\documentclass[12pt,twoside]{article}
\usepackage[english]{babel}
\usepackage[latin1]{inputenc}
\usepackage{amsmath}
\usepackage{amssymb,amsfonts}
\usepackage{graphicx}                   

\newcommand{\cT}{\mathcal{T}}

\newcommand{\HYP}{\mathbb{H}^3}
\newcommand{\HYN}{\mathbb{H}^n}

\begin{document}
\pagestyle{myheadings}
\markboth{\centerline{Jen\H o Szirmai}}
{Hyperball packings in hyperbolic $3$-space}
\title
{Hyperball packings in hyperbolic $3$-space}

\author{\normalsize{Jen\H o Szirmai} \\
\normalsize Budapest University of Technology and \\
\normalsize Economics Institute of Mathematics, \\
\normalsize Department of Geometry \\
\date{\normalsize{\today}}}

\maketitle


\begin{abstract}
In earlier works \cite{Sz06-1}, \cite{Sz06-2}, \cite{Sz13-3} and \cite{Sz13-4} we have investigated
the densest packings and the least dense coverings by congruent hyperballs (hyperspheres) to the regular prism
tilings in $n$-dimensional hyperbolic space $\HYN$ ($ 3 \le n \in \mathbb{N})$.

In this paper we study a large class of hyperball (hypersphere) packings in
$3$-dimensional hyperbolic space that can be derived from truncated simplex tilings (e.g. \cite{S14}, \cite{MPSz}).
It is clear, that in order to get a density upper bound for the above hyperball packings, it is sufficient to determine
the density upper bound locally, e.g. in truncated simplices.

Thus we study hyperball packings in truncated simplices, i.e. truncated tetrahedra and prove that if the truncated tetrahedron is regular, then the
density of the densest packing is $\approx 0.86338$. This is larger than the B\"or\"oczky-Florian density upper bound for balls and horoballs
(horospheres) \cite{B--F64} but our locally optimal hyperball packing configuration cannot be extended to the entirety of
hyperbolic space $\mathbb{H}^3$. But our regular truncated tetrahedron construction, under
the extended Coxeter group $[3, 3, 7]$ with top density $\approx 0.82251$, seems to be good enough (Table 1).

Moreover, we show that the densest known hyperball packing, related to the regular $p$-gonal prism tilings \cite{Sz06-1}, can dually be realized by
regular truncated tetrahedron
tilings as well.
\end{abstract}

\newtheorem{theorem}{Theorem}[section]
\newtheorem{corollary}[theorem]{Corollary}
\newtheorem{conjecture}{Conjecture}[section]
\newtheorem{lemma}[theorem]{Lemma}
\newtheorem{exmple}[theorem]{Example}
\newtheorem{defn}[theorem]{Definition}
\newtheorem{rmrk}[theorem]{Remark}
\newenvironment{definition}{\begin{defn}\normalfont}{\end{defn}}
\newenvironment{remark}{\begin{rmrk}\normalfont}{\end{rmrk}}
\newenvironment{example}{\begin{exmple}\normalfont}{\end{exmple}}
\newenvironment{acknowledgement}{Acknowledgement}


\section{Introduction}
Let $X$ denote either the $n$-dimensional sphere $\mathbb{S}^n$, Euclidean space $\mathbb{E}^n$, or hyperbolic space $\mathbb{H}^n$ $(2 \le n \in \mathbb{N})$
of constant curvature.

In space $X$ let $d_n(r)$ be the density of $n+1$ mutually touching spheres of radius $r$ with respect
to the simplex spanned by their centres. L.~Fejes T\'oth and H.~S.~M.~Coxeter
conjectured that the packing density of balls of radius $r$ in $X$ cannot exceed $d_n(r)$.
This conjecture has been proved by C.~A.~Rogers for Euclidean space $\mathbb{E}^n$ \cite{Ro64}.
The 2-dimensional spherical case was settled by L.~Fejes T\'oth in \cite{FTL}. In $\HYN$ there are $3$ kinds
of "generalized balls (spheres)":  the usual balls (spheres), horoballs (horospheres) and hyperballs (hyperspheres).

In \cite{B78} K.~B\"or\"oczky proved the following generalization
for {\it ball and horoball} packings for $n=3$, and claimed the analogous statement for any $n$:
\begin{theorem}[K.~B\"or\"oczky]
In an $n$-dimensional space of constant curvature consider a packing of spheres of radius $r$. In spherical space suppose that $r<\frac{\pi}{4}$.
Then the density of each sphere in its Dirichlet-Voronoi cell cannot exceed the density of $n+1$ spheres of radius $r$ mutually
touching one another with respect to the simplex spanned by their centers.
\end{theorem}

The above greatest density in $\mathbb{H}^3$ is $\approx 0.85328$
which is not realized by packing with equal balls. However, it is attained by the horoball packing of
$\overline{\mathbb{H}}^3$ where the ideal centers of horoballs lie on the
absolute figure of $\overline{\mathbb{H}}^3$. This ideal regular
simplex tiling is given with Coxeter-Schl\"afli symbol $[3,3,6]$.
Ball packings of hyperbolic $n$-space are extensively discussed in the literature see e.g. \cite{Be}, \cite{B78}, \cite{G--K--K},\cite{K98}
and  \cite{J}.

In a previous paper \cite{KSz} we proved that the above known optimal horoball packing arrangement in $\mathbb{H}^3$ is not unique.
We gave several new examples of horoball packing arrangements based on totally asymptotic Coxeter tilings that
yield the above B\"or\"oczky--Florian packing density upper bound \cite{B--F64}.
Furthermore, by admitting horoballs of different types at each vertex of a totally asymptotic simplex and generalizing
the simplicial density function to $\mathbb{H}^n$ $(n \ge 2)$, we have found that the B\"or\"oczky-Florian type density
upper bound is no longer valid for fully asymptotic simplices in higher dimensions $n > 3$  \cite{Sz12}, \cite{Sz12-2}.
For example in $\mathbb{H}^4$, the density of such optimal,
locally densest horoball packing is $\approx 0.77038$ larger than the
analogous B\"or\"oczky-Florian type density upper bound of $\approx 0.73046$.
However, these horoball packing configurations are only locally optimal and cannot be extended to the whole hyperbolic space $\mathbb{H}^4$.

In paper [14] we have continued our previous investigation in $\mathbb{H}^4$ allowing horoballs of different types.
We have shown seven counterexamples (which are realized by one-, two-, or three horoball types).
\newline
\indent In \cite{Sz06-1}, \cite{Sz06-2} and \cite{Sz13-3} we have studied the regular prism tilings and the corresponding optimal hyperball packings.
Their metric data and their densities have been determined. \newline
\indent In hyperbolic plane $\mathbb{H}^2$ the universal upper bound of the hypercycle packing density is $\frac{3}{\pi}$, and the universal lower bound of hypercycle covering density is $\frac{\sqrt{12}}{\pi}$,
proved by I.~Vermes in \cite{V79}. Recently, to the best of author's knowledge, candidates for the densest hyperball
(hypersphere) packings in the $3,4$ and $5$-dimensional hyperbolic spaces are derived by the regular prism
tilings \cite{Sz06-1}, \cite{Sz06-2} and \cite{Sz13-3}.

We observed that some extremal properties of hyperball packings naturally
belong to the regular truncated tetrahedron (or simplex, in general, see
Lemma 3.2 and Lemma 3.3).
Therefore, in this paper
we study hyperball packings in truncated tetrahedra, and prove that if the truncated tetrahedron is regular, then the
density of the densest packing is $\approx 0.86338$ (see Theorem 5.1). However, these hyperball packing configurations are only locally optimal,
and cannot be extended to the whole space $\mathbb{H}^3$.
Moreover, we
show that the densest known hyperball packing, dually related to the regular prism
tilings, introduced
in \cite{Sz06-1}, can be realized by a regular truncated tetrahedron tiling.

We have an extensive program of finding globally and locally optimal
ball packings in the eight Thurston geometries, arising from Thurston's geometrization conjecture
\cite{Sz07-1}, \cite{Sz07-2}, \cite{Sz10}, \cite{Sz12}, \cite{Sz12-2}, \cite{Sz13-1}, \cite{Sz13-2} and \cite{Sz14-1}.
Packing density is defined the as the ratio of the volume of a fundamental domain for
the symmetry group of a tiling to the volume of the ball pieces contained in
the fundamental domain. The large class of the truncated tetrahedron (or
simplex) tilings are studied, e.g. in \cite{S14}, on the base of \cite{MPSz}.

\section{The projective model and saturated hyperball packings in $\HYP$}
We use for $\mathbb{H}^3$ (and analogously for $\HYN$, $n>3$) the projective model in the Lorentz space $\mathbb{E}^{1,3}$
that denotes the real vector space $\mathbf{V}^{4}$ equipped with the bilinear
form of signature $(1,3)$,
$
\langle ~ \mathbf{x},~\mathbf{y} \rangle = -x^0y^0+x^1y^1+x^2y^2+ x^3 y^3,
$
where the non-zero vectors
$
\mathbf{x}=(x^0,x^1,x^2,x^3)\in\mathbf{V}^{4} \ \  \text{and} \ \ \mathbf{y}=(y^0,y^1,y^2,y^3)\in\mathbf{V}^{4},
$
are determined up to real factors, for representing points of $\mathcal{P}^n(\mathbb{R})$. Then $\mathbb{H}^3$ can be interpreted
as the interior of the quadric
$
Q=\{[\mathbf{x}]\in\mathcal{P}^3 | \langle ~ \mathbf{x},~\mathbf{x} \rangle =0 \}=:\partial \mathbb{H}^3
$
in the real projective space $\mathcal{P}^3(\mathbf{V}^{4},
\mbox{\boldmath$V$}\!_{4})$. Namely, for an interior point $\mathbf{y}$ holds $\langle ~ \mathbf{y},~\mathbf{y} \rangle <0$.

The points of the boundary $\partial \mathbb{H}^3 $ in $\mathcal{P}^3$
are called points at infinity, or at the absolute of $\mathbb{H}^3 $. The points lying outside $\partial \mathbb{H}^3 $
are said to be outer points of $\mathbb{H}^3 $ relative to $Q$. Let $P([\mathbf{x}]) \in \mathcal{P}^3$, a point
$[\mathbf{y}] \in \mathcal{P}^3$ is said to be conjugate to $[\mathbf{x}]$ relative to $Q$ if
$\langle ~ \mathbf{x},~\mathbf{y} \rangle =0$ holds. The set of all points which are conjugate to $P([\mathbf{x}])$
form a projective (polar) hyperplane
$
pol(P):=\{[\mathbf{y}]\in\mathcal{P}^3 | \langle ~ \mathbf{x},~\mathbf{y} \rangle =0 \}.
$
Thus the quadric $Q$ induces a bijection
(linear polarity $\mathbf{V}^{4} \rightarrow
\mbox{\boldmath$V$}\!_{4})$)
from the points of $\mathcal{P}^3$ onto their hyperplanes.

The point $X [\bold{x}]$ and the hyperplane $\alpha [\mbox{\boldmath$a$}]$
are incident if $\bold{x}\mbox{\boldmath$a$}=0$ ($\bold{x} \in \bold{V}^{4} \setminus \{\mathbf{0}\}, \ \mbox{\boldmath$a$}
\in \mbox{\boldmath$V$}_{4}
\setminus \{\mbox{\boldmath$0$}\}$).

The equidistance surface (or hypersphere) is a quadratic surface at a constant distance
from a plane (base plane) in both halfspaces. The infinite body of the hypersphere, containing the base plane, is called hyperball.

The {\it half hypersphere } with distance $h$ to a base plane $\pi$ is denoted by $\mathcal{H}^h_+$.
The volume of a bounded hyperball piece $\mathcal{H}^h_+(\mathcal{A})$,
delimited by a $2$-polygon $\mathcal{A} \subset \pi$, and its prism orthogonal to $\pi$, can be determined by the classical formula (2.1) of J.~Bolyai.
\begin{equation}
Vol(\mathcal{H}^h_+(\mathcal{A}))=\frac{1}{4}Area(\mathcal{A})\left[k \sinh \frac{2h}{k}+
2 h \right], \tag{2.1}
\end{equation}
The constant $k =\sqrt{\frac{-1}{K}}$ is the natural length unit in
$\mathbb{H}^3$. $K$ denotes the constant negative sectional curvature. In the following we may assume that $k=1$.

Let $\mathcal{H}^h$ be a hyperball packing in $\HYP$ with congruent hyperballs of height $h$. The density of any packing can be
improved by adding hyperballs as long as there is sufficient room to do so. Else we say that the packing is saturated. We always assume that
our packings are saturated. Sometimes we take a set of hyperballs $\{ \mathcal{H}^h_i\}$ of the hyperball packing $\mathcal{H}^h$.
Their base planes are denoted by $\beta_i$.
Thus in a saturated hyperball packing the distance between two ultraparallel base planes
$d(\beta_i,\beta_j)$ is at least $2h$ (where for the natural indices holds $i < j$
and $d$ is the hyperbolic distance function).
\section{On hyperball packings in a truncated tetrahedron}
We consider a saturated hyperball packing $\mathcal{H}^h$ of hyperballs in $\HYP$
which can be derived from a truncated tetrahedron tiling $\mathcal{T}$. One tetrahedron of $\mathcal{T}$ is
$\mathcal{S}$=$C_1^1 C_2^1 C_3^1$ $C_1^2 C_2^2 C_3^2$
$C_1^3 C_2^3 C_3^3$ $C_1^4 C_2^4 C_3^4$
illustrated in Fig.~1.~a.

The ultraparallel base planes of $\mathcal{H}^h_i$ $(i=1,2,3,4)$ are denoted by $\beta_i$. The distance between two base planes
$d(\beta_i,\beta_j)=:e_{ij} \ge 2h$ ($i < j \in \{1,2,3,4\})$.
Moreover, let the volume of the truncated tetrahedron (or simplex, in general) $\mathcal{S}$ be $Vol(\mathcal{S})$. We introduce the local density function
$\delta(\mathcal{S}(h))$ related to $\mathcal{S}$:
\begin{definition}
\begin{equation}
\delta(\mathcal{S}(h)):=\frac{\sum_{i=1}^4 Vol(\mathcal{H}^h_i \cap \mathcal{S})}{Vol({\mathcal{S}})}. \tag{3.1}
\end{equation}
\end{definition}
It is clear, that $\sup_{\mathcal{S}\in \cT}\delta(\mathcal{S}(h))$ provides an universal upper bound to the considered hyperball packing
$\mathcal{H}^h$ in space
$\HYP$.  The problem of determining $\sup_{\mathcal{S}}\delta(\mathcal{S})$ seems to be complicated in general, but we can formulate
some important assertions.
\begin{enumerate}
\item The area of each rectangular hexagon face of $\mathcal{S}$ is, e.g. $$Area(C_1^1 C_2^1 C_1^3 C_2^3 C_2^2 C_1^2)=\pi.$$
\item If we restrict ourselves to the above rectangular hexagon $\mathcal{F}=C_1^1 C_2^1 $ $C_1^3 C_2^3 C_2^2 C_1^2$ then
the intersections of $\mathcal{H}^h_i$ $(i=1,2,3)$ with $\mathcal{F}$ form in $\mathcal{F}$ a partial hypercycle packing (see Fig.~1.~b).

It is clear, that the density $\delta(\mathcal{F}(h))$ of the hypercycle packing in $\mathcal{F}$ is maximal if the area
$\sum_{i=1}^3 Area(\mathcal{H}^h_i \cap \mathcal{F})$ is maximal, because $Area(\mathcal{F})=\pi$ is fixed.
I.~Vermes in \cite{V79} noticed that the density $\delta(\mathcal{F}(h))$ is maximal if the lengths of the common perpendiculars are equal to
$e_{12}=e_{23}=e_{13}=2h$. We note here, that in this "regular" case $\sum_{i=1}^3(b_i)$ is maximal, as well, where $b_i$ are the "base segments" of
the hypercycle domains $\mathcal{H}^h_i \cap \mathcal{F}$.
I. Vermes proved in \cite{V79} that
$$\delta(\mathcal{F}(h))=\frac{6\sinh({h})\mathrm{arcsinh}\frac{1}{2\sinh({h})}}{\pi},
~~~~~\lim_{h\rightarrow\infty}(\delta(\mathcal{F}(h)))=\frac{3}{\pi},~ \text{increasingly.}$$
\end{enumerate}
\begin{lemma}We can extend the above statement to the other rectangular hexagon facets. Therefore, if the distance between two base planes is
$e_{ij}=2h$ ($i < j \in \{1,2,3,4\})$, then the "regular" truncated tetrahedron provides the densest hypercycle packing
in the rectangular hexagons of $\mathcal{S}$, and the density of the densest hypercycle packing can be improved in the rectangular hexagon
facets with the density of the above hypercycle packings if $h \rightarrow \infty$.
\end{lemma}
The dihedral angles of the truncated tetrahedron $\mathcal{S}$ at the edges $B_iB_j$, ($i<j \in \{1,2,3,4\}$ are denoted by $\omega_{ij}$.
In the following assume that the sum of the dihedral angles $\omega_{ij}$ is constant: $\Omega$. (At truncations the
other dihedral angles of $\mathcal{S}$ are $\frac{\pi}{2}$). We obtain the following
lemma as a consequence of the above assertions and formula (2.1).
\begin{lemma}
If the sum of the above dihedral angles $\omega_{ij}$ is constant: $\Omega$, then the surface area of
$\mathcal{S}$ is $8\pi-2\Omega$ constant as well. Therefore,
$\sum_{k=1}^4 Vol(\mathcal{H}^h_k \cap \mathcal{S})$ is maximal if $e_{ij}=2h$ $(i<j \in\{1,2,3,4 \}$.
\end{lemma}
Although this lemma does not provide explicite estimate yet ($h$ depends on $\omega_ij$'s),
it motivates the following additional assumption: let the truncated tetrahedron be regular.
Then $h$ can also be calculated as follows.
\begin{figure}[ht]
\centering
\includegraphics[width=6cm]{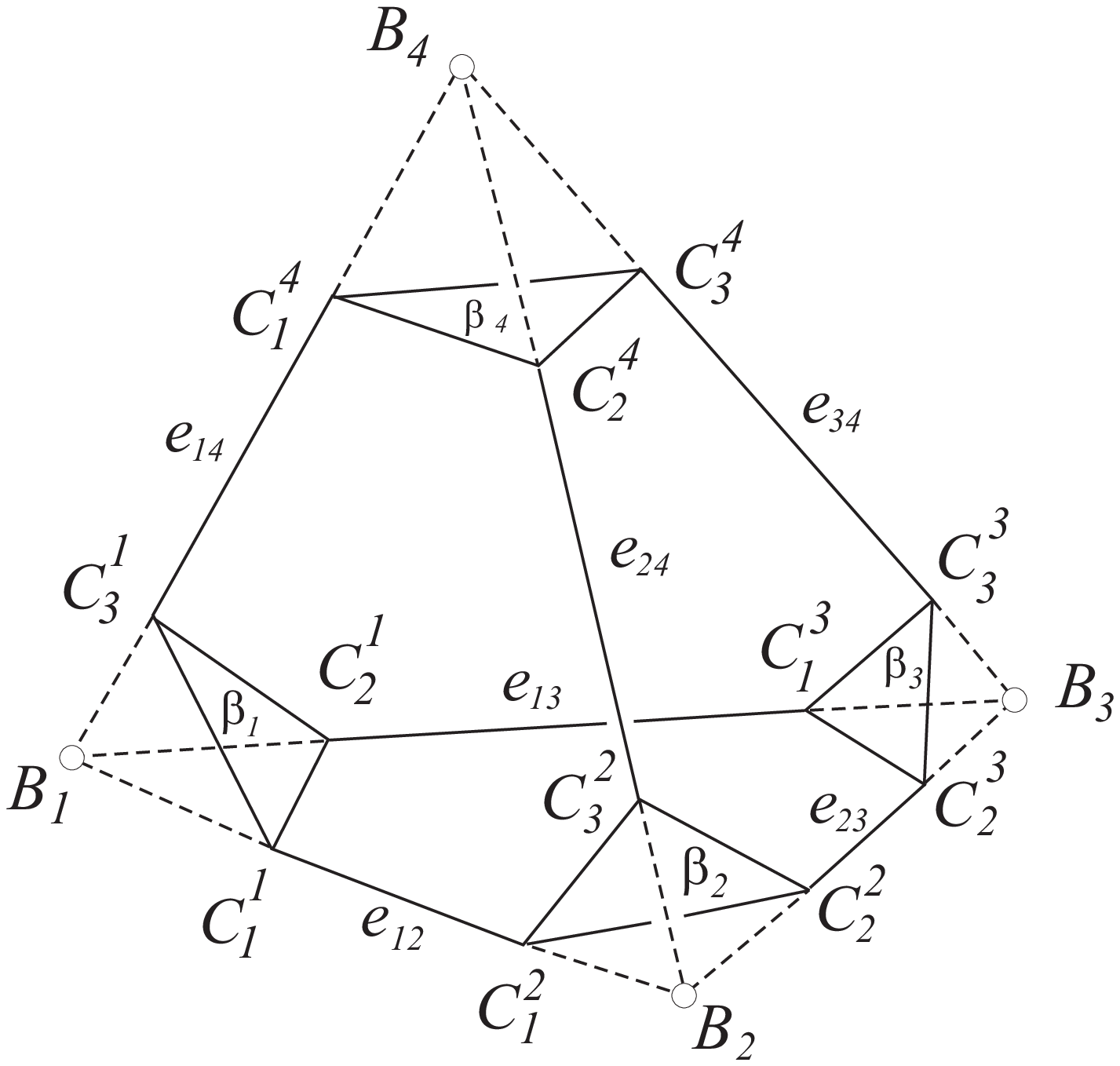} \includegraphics[width=6cm]{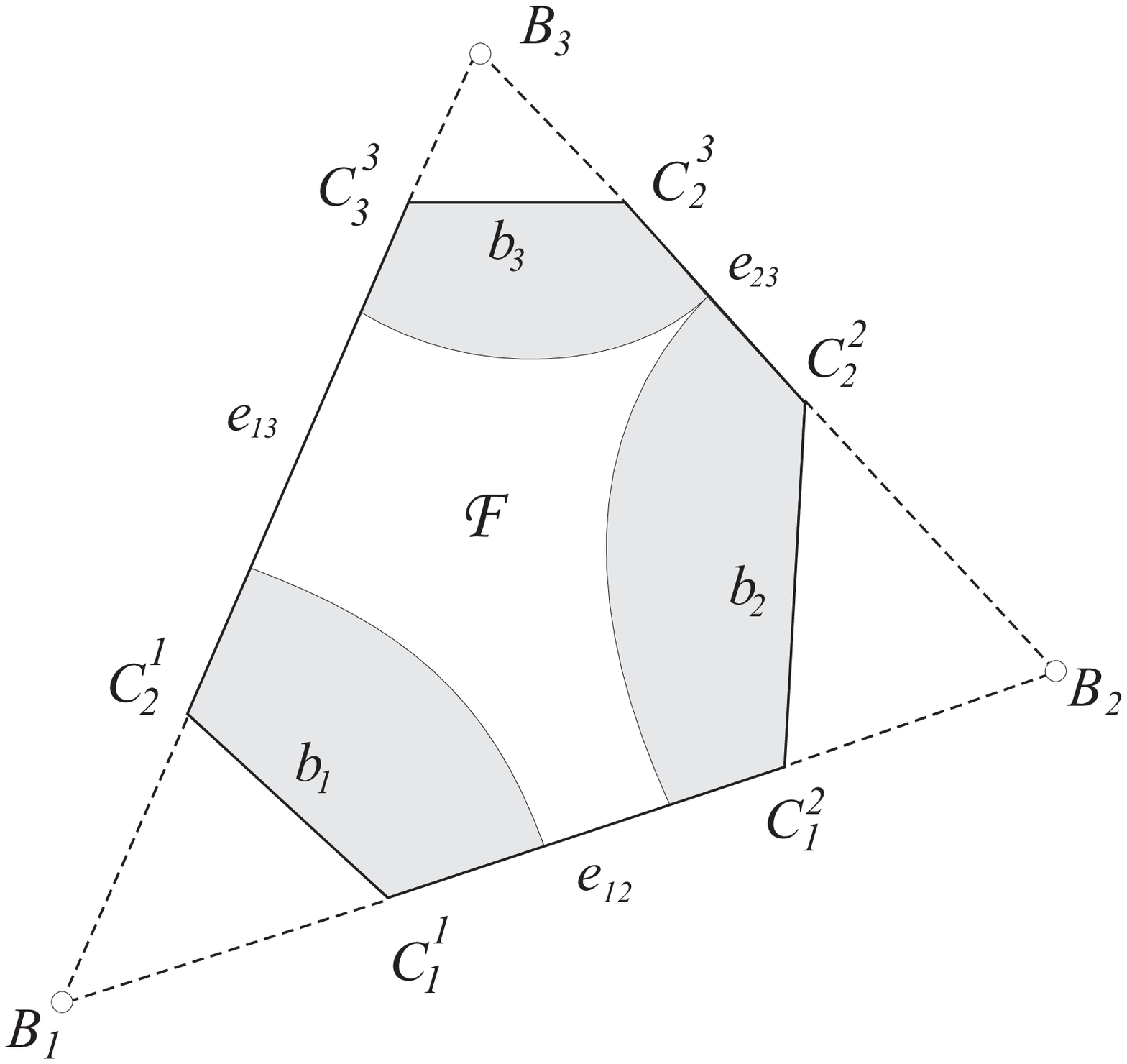}

a.~~~~~~~~~~~~~~~~~~~~~~~~~~~~~~~~~~~~b.
\caption{Truncated tetrahedron and one of its rectangular hexagon faces}
\label{}
\end{figure}
\section{Characteristic orthoschemes for the volume of a truncated regular tetrahedron}

{\it An orthoscheme $\mathcal{O}$ in $\mathbb{H}^n$ $n \geq 2$ in classical sense}
is a simplex bounded by $n+1$ hyperplanes $H_0,\dots,H_n$
such that (\cite{B--H}) $H_i \bot H_j, \  \text{for} \ j\ne i-1,i,i+1.$
\begin{rmrk}
This definition is equivalent to the following: A simplex $\mathcal{O}$ in $\mathbb{H}^n$ is an
orthoscheme iff the $n+1$ vertices of $\mathcal{O}$ can be
labelled by $R_0,R_1,\dots,R_n$ in such a way that
$\text{span}(R_0,\dots,R_i) \perp \text{span}(R_i,\dots,R_n) \ \ \text{for} \ \ 0<i<n-1.$
\end{rmrk}

Geometrically, complete orthoschemes of degree $m$ can be described as
follows:

\begin{enumerate}
\item
For $m=0$, they coincide with the class of classical orthoschemes introduced by
{{Schl\"afli}}. The initial and final vertices, $R_0$ and $R_n$ of the orthogonal edge-path
$R_iR_{i+1},~ i=0,\dots,n-1$, are called principal vertices of the orthoscheme (see Remark 4.1).
\item
A complete orthoscheme of degree $m=1$ can be interpreted as an
orthoscheme with one outer principal vertex, say $R_n$, which is truncated by
its polar plane $pol(R_n)$ (see Fig.~2.~b). In this case the orthoscheme is called simply truncated with
outer vertex $R_n$.
\item
A complete orthoscheme of degree $m=2$ can be interpreted as an
orthoscheme with two outer principal vertex, $R_0,~R_n$, which is truncated by
its polar hyperplanes $pol(R_0)$ and $pol(R_n)$. In this case the orthoscheme is called doubly
truncated (see \cite{K89}).
\end{enumerate}
\begin{figure}[ht]
\centering
\includegraphics[width=5.5cm]{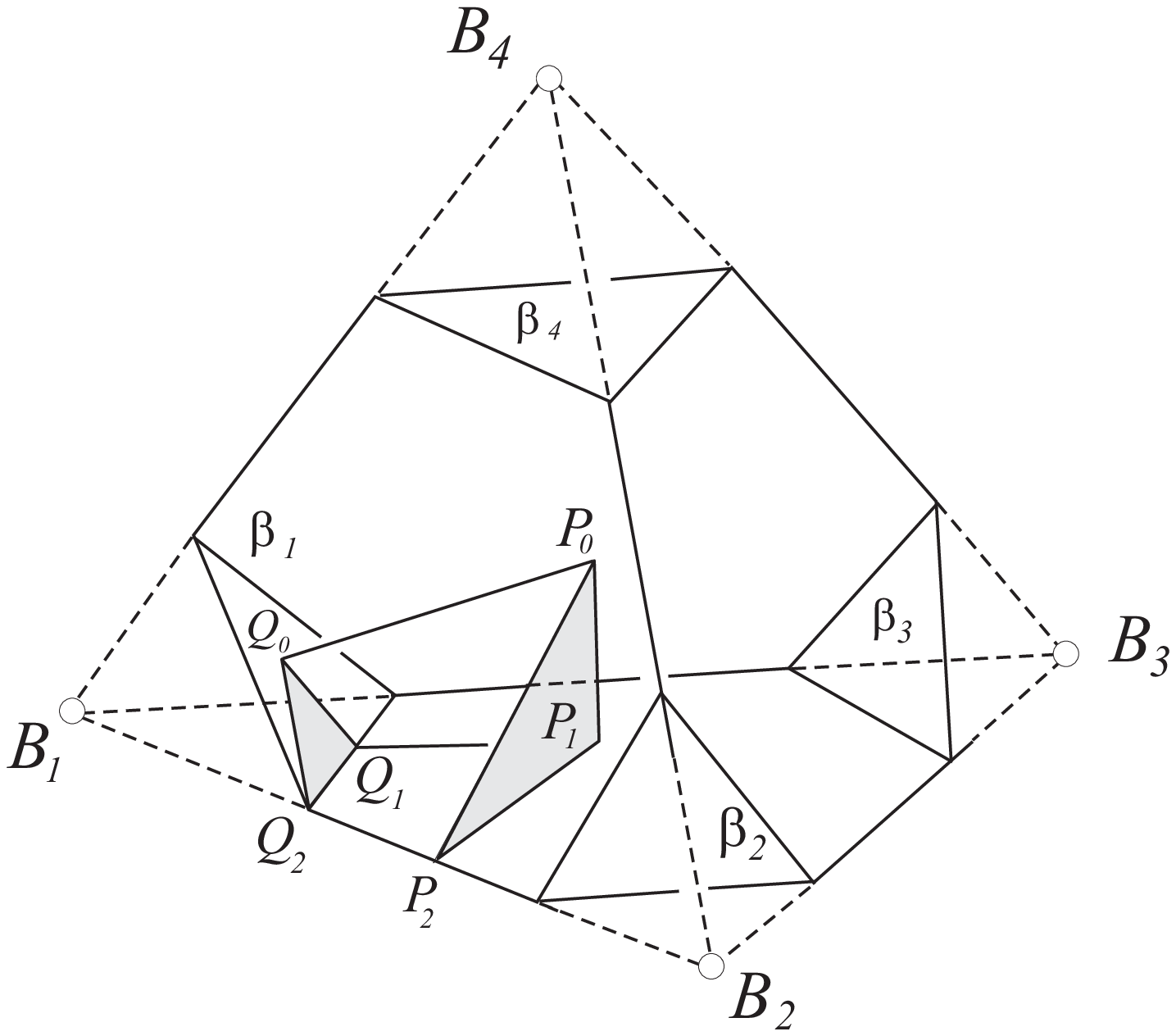} \includegraphics[width=6.5cm]{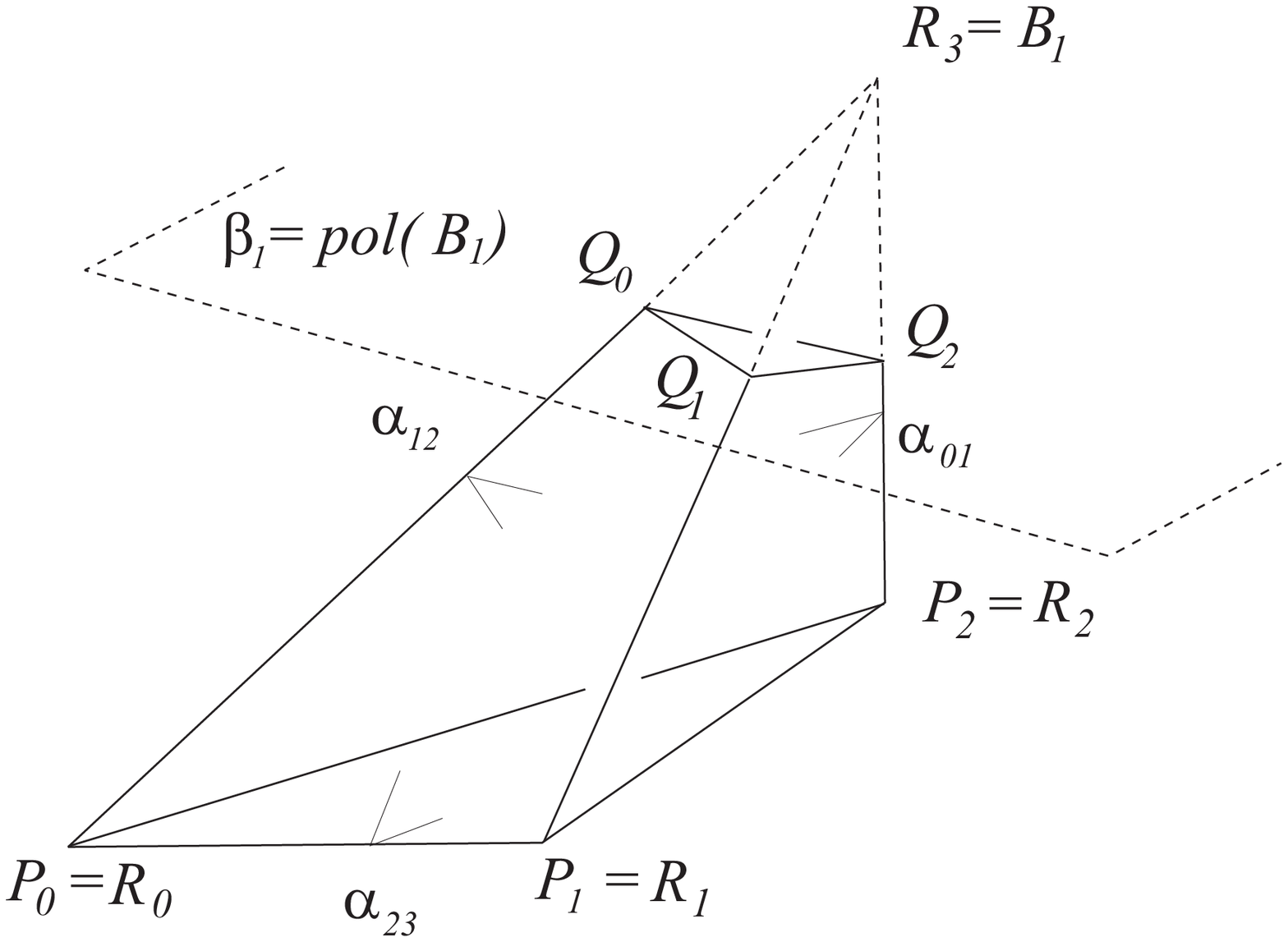}

a.~~~~~~~~~~~~~~~~~~~~~~~~~~~~~~~~~~~~~b.
\caption{Truncated tetrahedron with a complete orthoscheme of degree $m=1$ (simple frustum orthoscheme)}
\label{}
\end{figure}
In the following we use the "$3$-dimensional simple frustum orthoschemes" whose volume formula is derived by the next
Theorem of {{R.~Kellerhals}} \cite{K89}:
\begin{theorem}{\rm{(R.~Kellerhals)}} The volume of a three-dimensional hyperbolic
complete orthoscheme $\mathcal{O} \subset \mathbb{H}^3$
is expressed with the essential
angles $\alpha_{01},\alpha_{12},\alpha_{23}, \ (0 \le \alpha_{ij} \le \frac{\pi}{2})$
(Fig.~1.~b) in the following form:

\begin{align}
&Vol(\mathcal{O})=\frac{1}{4} \{ \mathcal{L}(\alpha_{01}+\theta)-
\mathcal{L}(\alpha_{01}-\theta)+\mathcal{L}(\frac{\pi}{2}+\alpha_{12}-\theta)+ \notag \\
&+\mathcal{L}(\frac{\pi}{2}-\alpha_{12}-\theta)+\mathcal{L}(\alpha_{23}+\theta)-
\mathcal{L}(\alpha_{23}-\theta)+2\mathcal{L}(\frac{\pi}{2}-\theta) \}, \notag
\end{align}
where $\theta \in [0,\frac{\pi}{2})$ is defined by:
$$
\tan(\theta)=\frac{\sqrt{ \cos^2{\alpha_{12}}-\sin^2{\alpha_{01}} \sin^2{\alpha_{23}
}}} {\cos{\alpha_{01}}\cos{\alpha_{23}}},
$$
and where $\mathcal{L}(x):=-\int\limits_0^x \log \vert {2\sin{t}} \vert dt$ \ denotes the
Lobachevsky function.
\end{theorem}

In the following we assume that the ultraparallel base planes $\beta_i$ of $\mathcal{H}^{h(p)}_i$ $(i=1,2,3,4)$
generate a "regular truncated tetrahedron" $\mathcal{S}^r$ with outer vertices $B_i$ (see Fig.~2.~a)
the non-orthogonal dihedral angles of $\mathcal{S}^r$ are equal to $\frac{\pi}{p}$,
$(6<p\in \mathbb{R})$ and the distances between two base planes $d(\beta_i,\beta_j)=:e_{ij}$  ($i < j \in \{1,2,3,4\})$ are equal to $2h(p)$ depending on the angle
$\frac{\pi}{p}$.

The truncated regular tetrahedron $\mathcal{S}^r$ can be decomposed into $24$ congruent simple frustum orthoschemes; one of them
$\mathcal{O}=Q_0Q_1Q_2P_0P_1P_2$ is illustrated in Fig.~2.~a where $P_0$ is the centre of the "regular tetrahedron" $\mathcal{S}^r$,
$P_1$ is the centre of a hexagonal face of $\mathcal{S}^r$, $P_0$ is the midpoint of a "common perpendicular" edge of this face,
$Q_0$ is the centre of an adjacent regular triangle face of $\mathcal{S}^r$, $Q_1$ is the midpoint of an appropriate edge of this face and
one of its endpoint is $Q_2$.

In our case the dihedral angles of orthoschemes $\mathcal{O}$ are the following:
$\alpha_{01}=\frac{\pi}{p}, \ \ \alpha_{12}=\frac{\pi}{3}, \ \
\alpha_{23}=\frac{\pi}{3}$ (see Fig.~2.~b). Therefore, the volume $Vol(\mathcal{O})$ of the orthoscheme $\mathcal{O}$ and the volume
$Vol(\mathcal{S}^r)=24 \cdot Vol(\mathcal{O})$ can be computed for any given parameter $p$ $(6<p\in \mathbb{R})$ by Theorem 4.2.
\section{Packing with congruent hyperballs in a regular truncated tetrahedron}
In this case for a given parameter $p$ the length of the common perpendiculars $h(p)=\frac{1}{2}e_{ij}$ $(i < j$, $i,j \in \{1,2,3,4\})$
can be determined by the machinery of projective metric geometry.

The points $P_2[{\mathbf{p}}_2]$ and $Q_2[{\mathbf{q}}_2]$ are proper points of hyperbolic $3$-space and
$Q_2$ lies on the polar hyperplane $pol(B_1)[\mbox{\boldmath$b$}^1]$ of the outer point $B_1$ thus
\begin{equation}
\begin{gathered}
\mathbf{q}_2 \sim c \cdot \mathbf{b}_1 + \mathbf{p}_2 \in \mbox{\boldmath$b$}^1 \Leftrightarrow
c \cdot \mathbf{b}_1 \mbox{\boldmath$b$}^1+\mathbf{p}_2 \mbox{\boldmath$b$}^1=0 \Leftrightarrow
c=-\frac{\mathbf{p}_2 \mbox{\boldmath$b$}^1}{\mathbf{b}_1 \mbox{\boldmath$b$}^1} \Leftrightarrow \\
\mathbf{q}_2 \sim -\frac{\mathbf{p}_2 \mbox{\boldmath$b$}^1}{\mathbf{b}_1 \mbox{\boldmath$b$}^1}
\mathbf{b}_1+\mathbf{p}_2 \sim \mathbf{p}_2 (\mathbf{b}_1 \mbox{\boldmath$b$}^1) - \mathbf{b}_1 (\mathbf{p}_2 \mbox{\boldmath$b$}^1)=
\mathbf{p}_2 h_{33}-\mathbf{b}_1 h_{23},
\end{gathered} \tag{5.1}
\end{equation}
where $h_{ij}$ is the inverse of the Coxeter-Schl\"afli matrix
\[(c^{ij}):=\begin{pmatrix}
1& -\cos{\frac{\pi}{p}}& 0&0\\
-\cos{\frac{\pi}{p}} & 1 & -\cos{\frac{\pi}{3}}&0\\
0 & -\cos{\frac{\pi}{3}} & 1&-\cos{\frac{\pi}{3}} \\
0&0&-\cos{\frac{\pi}{3}}&1\\
\end{pmatrix} \tag{5.2}
\]
of the orthoscheme $\mathcal{O}$.
The hyperbolic distance $h(p)$ can be calculated by the following formula:
\[
\begin{gathered}
\cosh{h(p)}=\cosh{P_2Q_2}=\frac{- \langle {\mathbf{q}}_2, {\mathbf{p}}_2 \rangle }
{\sqrt{\langle {\mathbf{q}}_2, {\mathbf{q}}_2 \rangle \langle {\mathbf{p}}_2, {\mathbf{p}}_2 \rangle}}= \\ =\frac{h_{23}^2-h_{22}h_{33}}
{\sqrt{h_{22}\langle \mathbf{q}_2, \mathbf{q}_2 \rangle}} =
\sqrt{\frac{h_{22}~h_{33}-h_{23}^2}
{h_{22}~h_{33}}}.
\end{gathered} \notag
\]
We get that the volume $Vol(\mathcal{S}^r)$, the maximal height $h(p)$ of the congruent hyperballs lying in $\mathcal{S}^r$ and
$\sum_{i=1}^4 Vol(\mathcal{H}^h_i \cap \mathcal{S}^r))$ depends only on the parameter $p$ of the truncated regular tetrahedron $\mathcal{S}^r$.

Therefore, the density $\delta(\mathcal{S}^r(h(p)))$ depends only on $p$ $(6<p\in \mathbb{R})$. Moreover,
the volume of the hyperball pieces can be computed by formula (2.1), and the volume of $\mathcal{S}^r$ can be determined by Theorem 4.2.
\begin{figure}[ht]
\centering
\includegraphics[width=6cm]{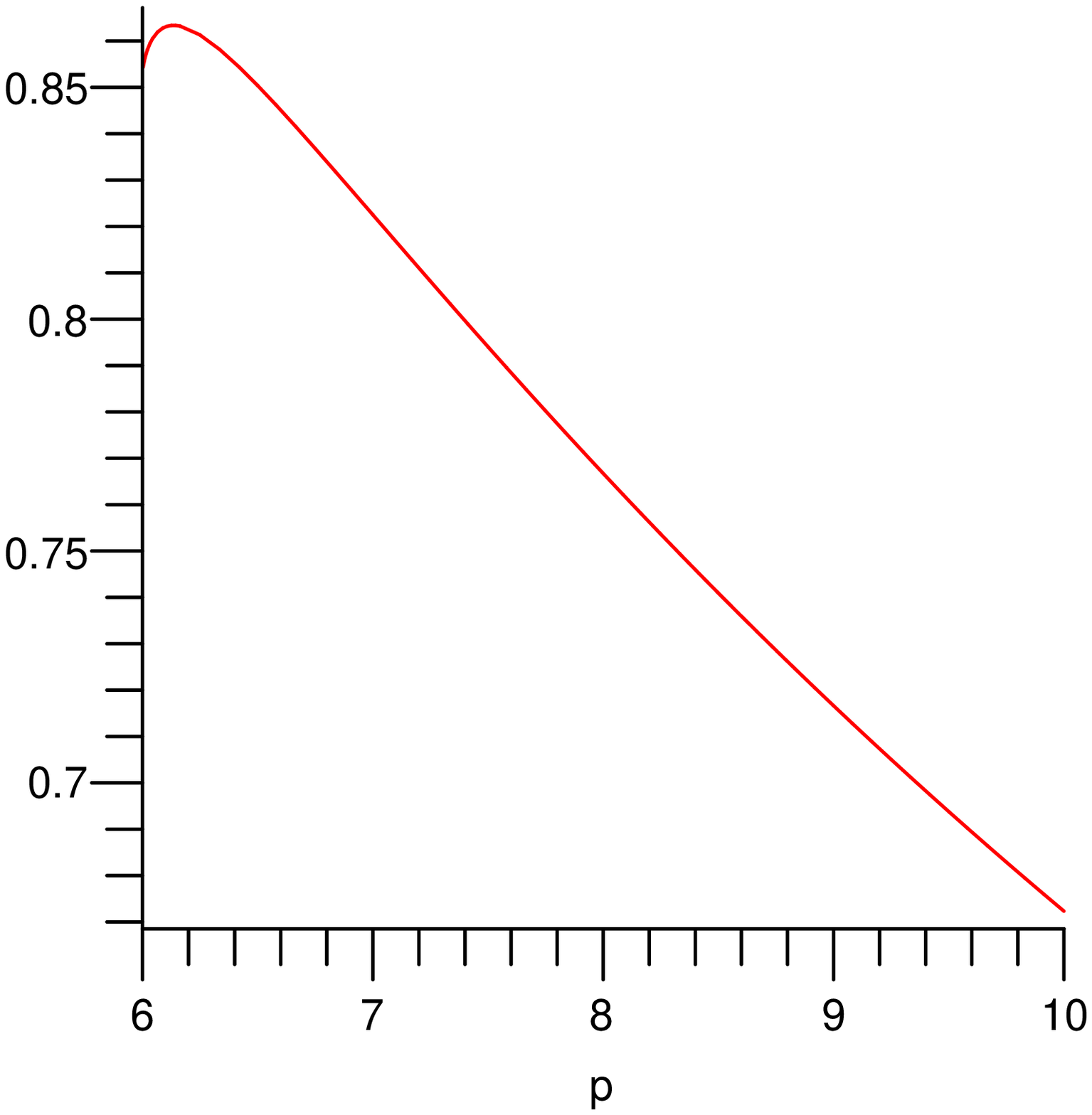} \includegraphics[width=6cm]{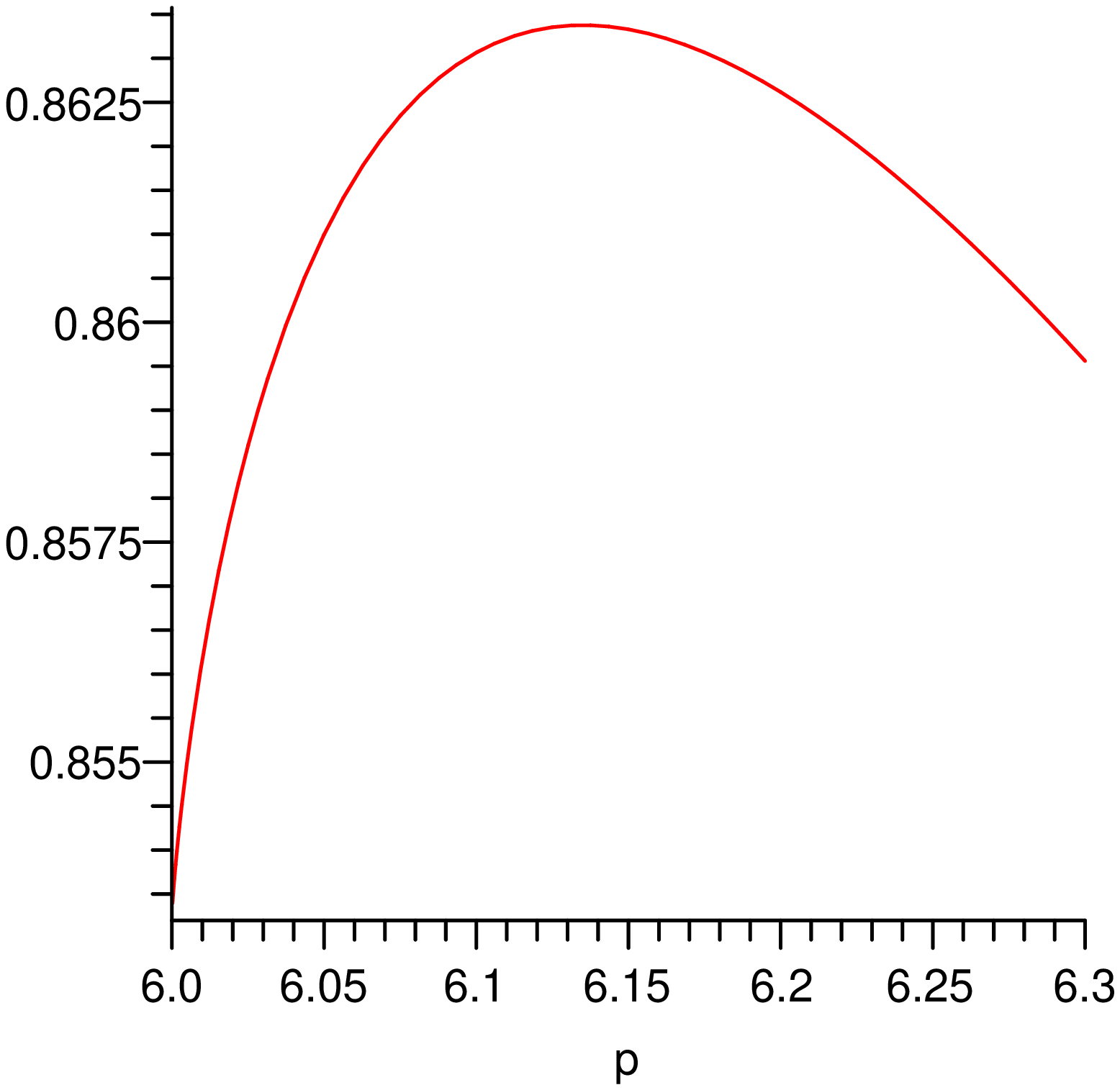}
\caption{}
\label{}
\end{figure}
Finally, we obtain after careful analysis of the smooth density function the following
\begin{theorem}
The density function $\delta(\mathcal{S}^r(h(p)))$, $p\in (6,\infty)$
attains its maximum at $p^{opt} \approx 6.13499$, and $\delta(\mathcal{S}^r(h(p)))$
is strictly increasing in the interval $(6,p^{opt})$, and strictly decreasing in $(p^{opt},\infty)$. Moreover, the optimal density
$\delta^{opt}(\mathcal{S}^r(h(p^{opt}))) \approx 0.86338$ (see Fig.~3).
\end{theorem}
\begin{rmrk}
\begin{enumerate}
\item In our case $\lim_{p\rightarrow 6}(\delta(\mathcal{S}^r(h(p))))$ is equal to the B\"oröczky-Florian
upper bound of the ball and horoball packings in $\HYP$ \cite{B--F64}.
\item $\delta^{opt}(\mathcal{S}^r(h(p^{opt}))) \approx 0.86338$ is larger than the B\"oröczky-Florian upper bound;
$\approx 0.85328$; but these hyperball packing configurations
are only locally optimal and cannot be extended to the entirety of hyperbolic space $\mathbb{H}^3$.
\end{enumerate}
\end{rmrk}
\subsection{Tilings with regular truncated tetrahedra in hyperbolic $3$-space}
In papers  \cite{Sz13-1}, \cite{Sz13-2}, \cite{Sz06-1}, \cite{Sz06-2} we have studied the hyperball packings and coverings
to regular prism tilings in $n$-dimensional $(n=3,4,5)$ hyperbolic space
and determined the corresponding densest hyperball packings and least dense hyperball coverings.
From the definitions of the prism tilings and the complete orthoschemes of degree $m=1$ it follows that a
regular prism tiling exists in space $\mathbb{H}^n$ if and only if there
exists a complete Coxeter orthoscheme of degree $m=1$ with two ultraparallel faces.
The complete Coxeter orthoschemes were classified by {{Im Hof}} in
\cite{IH85} and \cite{IH90} by generalizing the methods of {{Coxeter}} and {{B\"ohm}} appropriately.
The truncated tetrahedron tilings are studied e.g. in \cite{S14} on the base of \cite{MPSz}.

The hyperball packings in the regular truncated tetrahedra under the extended reflection groups with Coxeter-Schläfli symbol
$[3,3,p]$, investigated in this paper, can be extended to the entire hyperbolic space
if $6<p$ integer parameter and coincide with the hyperball packings to the regular $p$-gonal prism tilings in $\HYP$ with extended
Coxeter-Schl\"afli symbols
$[p,3,3]$, which are discussed in \cite{Sz06-1}, because their vertex figure is the tetrahedron given by Schl\"afli symbol [3,3].
As we know, $[3,3,p]$ and $[p,3,3]$ are dually isomorphic extended reflection groups, just with the above frustum orthoscheme
as fundamental domain (Fig.~2.b, matrix $(c^{ij})$ in formula (5.2)).

In the following Table we summarize the data of the hyperball packings for some parameters $p$, ($6<p\in \mathbb{N}$).
\medbreak
\centerline{\vbox{
\halign{\strut\vrule~\hfil $#$ \hfil~\vrule
&\quad \hfil $#$ \hfil~\vrule
&\quad \hfil $#$ \hfil\quad\vrule
&\quad \hfil $#$ \hfil\quad\vrule
&\quad \hfil $#$ \hfil\quad\vrule
\cr
\noalign{\hrule}
\multispan5{\strut\vrule\hfill\bf Table 1, \hfill\vrule}%
\cr
\noalign{\hrule}
\noalign{\vskip2pt}
\noalign{\hrule}
p & h & Vol(\mathcal{O}) & Vol(\mathcal{H}^h_+(\mathcal{A}))& \delta^{opt} \cr
\noalign{\hrule}
7 & 0.78871 & 0.08856 & 0.07284 & 0.82251 \cr
\noalign{\hrule}
8 & 0.56419 & 0.10721 & 0.08220 & 0.76673 \cr
\noalign{\hrule}
9 & 0.45320 & 0.11825 & 0.08474 & 0.71663 \cr
\noalign{\hrule}
\vdots & \vdots  & \vdots  & \vdots  & \vdots \cr
\noalign{\hrule}
20 & 0.16397 & 0.14636 & 0.06064 & 0.41431 \cr
\noalign{\hrule}
\vdots & \vdots  & \vdots  & \vdots  & \vdots \cr
\noalign{\hrule}
50 & 0.06325 & 0.15167 & 0.02918 & 0.19240 \cr
\noalign{\hrule}
\vdots & \vdots  & \vdots  & \vdots  & \vdots \cr
\noalign{\hrule}
100 & 0.03147 & 0.15241 & 0.01549 & 0.10165 \cr
\noalign{\hrule}
p \to \infty & 0 & 0.15266 & 0 & 0 \cr
\noalign{\hrule}}}}
\medbreak
The problems of the densest horoball and hyperball packings in hyperbolic $n$-space $n \ge 3$ with horoballs
of different types and hyperballs has not been settled yet, in general (see e.g. \cite{KSz}, \cite{KSz14}, \cite{Sz12}, \cite{Sz12-2}).

Optimal ball (sphere) packings in other homogeneous Thurston geometries represent
another huge class of open problems. For these non-Euclidean geometries
only very few results are known (e.g. \cite{Sz07-2}, \cite{Sz10}, \cite{Sz13-2}, \cite{Sz14-1}). Detailed studies are the objective of ongoing research.
Applications of the above projective method seem to be interesting in (non-Euclidean)
crystallography as well, a topic of much current interest.

I thank Prof. Emil Moln\'ar for helpful comments and improvements to this paper.

\noindent
\footnotesize{Budapest University of Technology and Economics Institute of Mathematics, \\
Department of Geometry, \\
H-1521 Budapest, Hungary. \\
E-mail:~szirmai@math.bme.hu \\
http://www.math.bme.hu/ $^\sim$szirmai}

\end{document}